\documentclass[proceedings,submission%
]{dmtcs}

\usepackage[latin1]{inputenc}

\newlength{\cellsize}
\title{A $t$-generalization for Schubert Representatives of the Affine Grassmannian}

\author{Avinash J. Dalal
\thanks{Partially supported by the NSF grant DMS--1001898.}
and 
Jennifer Morse
\thanks{Partially supported by the NSF grants DMS--1001898.}}

\address{Department of Mathematics, Drexel University, Philadelphia, PA 19104, U.S.A.}


\usepackage{amsmath,amssymb,latexsym,graphics,color,graphicx}
\usepackage[all]{xy}
\usepackage{tabmac_avi}
\usepackage[usenames,dvipsnames]{xcolor}

\newtheorem{theorem}{Theorem}
\newtheorem{property}[theorem]{Property}
\newtheorem{cor}[theorem]{Corollary}
\newtheorem{lemma}[theorem]{Lemma}

\newtheorem{defn}[theorem]{Definition}

\newtheorem{remark}[theorem]{Remark}
\newtheorem{exmpl}[theorem]{Example}
\newtheorem{example}[theorem]{Example}
\newtheorem{conj}[theorem]{Conjecture}

\def\core {\mathfrak c}

\def \core {{\mathfrak c}}

\def\Gr {\text{Gr}}

\def\charge{ {\rm {charge}}}

\def\weight{ {\rm {weight}}}

\def\shape{ {\rm {shape}}}

\begin{document}
\newcommand{\apnd}{A}
\newcommand{\rbar}{\bar{R_1}}
\newcommand{\boldx}{\textbf{x}}
\newcommand{\HRule}{\rule{\linewidth}{0.5mm}}
\newcommand{\emptbx}{\mbox{}}
\newcommand{\filbx}{\rule[-0.45mm]{4mm}{4mm}}

\maketitle
\begin{abstract}
\paragraph{Abstract.}
We introduce two families of symmetric functions with an extra parameter
$t$ that specialize to Schubert representatives for cohomology
and homology of the affine Grassmannian when $t=1$.  The
families are defined by a statistic on combinatorial objects associated
to the type-$A$ affine Weyl group and their transition matrix with
Hall-Littlewood polynomials is $t$-positive.  We conjecture that one
family is the set of $k$-atoms.

Nous pr\'esentons deux familles de fonctions sym\'etriques d\'ependant d'un
param\`etre $t$ et dont les sp\'ecialisations \`a $t=1$ correspondent 
aux classes de Schubert dans la cohomologie et l'homologie des vari\'et\'es Grassmanniennes affines.  Les familles sont d\'efinies par des statistiques sur certains objets combinatoires associ\'es au groupe de Weyl affine de type $A$ et leurs matrices de transition dans la base des polyn\^omes de Hall-Littlewood sont $t$-positives.  Nons conjecturons qu'une de ces familles 
correspond aux $k$-atomes.

\end{abstract}

\section{Introduction}
\label{sec:in}

{Affine Schubert calculus} is a generalization of classical
Schubert calculus where the Grassmannian is replaced by 
infinite-dimensional spaces $\Gr_G$ known as affine Grassmannians.  
As with Schubert calculus, topics under the umbrella of affine Schubert 
calculus are vast but now, it is the combinatorics of a family of polynomials 
called $k$-Schur functions that underpins the theory.

The theory of $k$-Schur functions came out of a
study of 
symmetric functions over $\mathbb Q(q,t)$ called
Macdonald polynomials.  Macdonald polynomials posses 
remarkable properties whose proofs inspired deep 
work in many areas 
One aspect that has been intensely studied
from a combinatorial, representation theoretic, and algebraic geometric 
perspective is the Macdonald/Schur transition matrix.
In particular,
in the late 1980's, Macdonald conjectured \cite{[M2]} that the coefficients
in the expansion
\begin{equation}
\label{macdo}
H_\mu[X;q,t] = \sum_\lambda K_{\lambda,\mu}(q,t) \,s_\lambda
\end{equation}
are positive sums of monomials
in $q$ and $t$; that is, $K_{\lambda,\mu}(q,t)\in\mathbb N[q,t]$.
These coefficients have since been a matter of great interest.
For starters, they generalize the {\it Kostka-Foulkes 
polynomials}.  These are given by  $K_{\lambda,\mu}(0,t)$
and they appear in many contexts such as Hall-Littlewood polynomials
\cite{Green}, affine Kazhdan-Lusztig theory \cite{Lu}, and affine tensor
product multiplicities~\cite{NY:1997}.  Moreover, Kostka-Foulkes 
polynomials encode the dimensions of certain bigraded $S_n$-modules \cite{GP}.  
They were beautifully characterized 
by Lascoux and Sch\"utzenberger \cite{LSfoulkes} 
by associating a {statistic} (non-negative integer) called {\it charge} 
to each tableau $T$ so that 
\begin{equation}
\label{charge}
K_{\lambda,\mu}(0,t)=
\sum_{{\substack{weight(T)= \mu \\\ \shape(T)=\lambda}}} t^{\charge(T)}
\,.
\end{equation}
Despite having such concrete 
results for the $q=0$ case, it was a big effort even to establish 
polynomiality for general $K_{\lambda,\mu}(q,t)$ 
\cite{GR,GT,Knop,LV,KiN,S} and the geometry of Hilbert 
schemes was eventually needed to prove positivity \cite{Haiman}.
A formula in the spirit of \eqref{charge} still remains a mystery.

In one study of Macdonald polynomials, Lapointe, Lascoux, and 
Morse found computational evidence for a family of new bases 
\begin{equation}
\label{katoms}
\{A_\mu^{(k)}[X;t]\}_{\mu_1\leq k}
\end{equation}
for subspaces 
$$
\Lambda_t^{(k)}= \text{span}\{H_\lambda[X;q,t]\}_{\lambda_1\leq k}\
$$
in a filtration $\Lambda_t^{(1)}\subseteq\Lambda_t^{(2)}\subseteq\dots
\subseteq\Lambda_t^{(\infty)}$ of  $\Lambda$.  Conjecturally, the
star feature of each basis was the property that Macdonald polynomials 
expand positively in terms of it, giving a remarkable factorization 
for the Macdonald/Schur transition matrices over $\mathbb N[q,t]$.  
To be precise, for any fixed integer $k>0$
and each $\lambda\in \mathcal P^k$ (a partition
where $\lambda_1\leq k$), 
\begin{equation}
H_{\lambda}[X;q,t\, ] = \sum_{\mu\in\mathcal P^k}
K_{\mu,\lambda}^{(k)}(q,t) \, A_{\mu}^{(k)}[X;t\, ] \quad\text{where}\quad
K_{\mu,\lambda}^{(k)}(q,t) \in \mathbb N[q,t] \, .
\label{mackkostka}
\end{equation}
It was conjectured in \cite{[LLM]} that 
for all $k>0$, $\{A^{(k)}_\mu[X;t]\}_{\mu_1\leq k}$ exists
and forms a basis for $\Lambda^{(k)}_t$,
and that for $k\geq |\mu|$, $A_\mu^{(k)}[X;t]=s_\mu$.
These conjectures and the decomposition \eqref{mackkostka} strengthen 
Macdonald's conjecture.  

A construction for $A_\mu^{(k)}[X;t]$ is given in 
\cite{[LLM]}, but it is so intricate that these conjectures 
remain unproven.  However,
pursuant investigations of these bases led to
various conjecturally equivalent characterizations.
One such family of polynomials $\{s_\lambda^{(k)}\}$
was introduced in \cite{[LMcore]} and conjectured to be the $t=1$ 
case of $A_\lambda^{(k)}[X;t]$.  It has since been proven that 
the $s_\lambda^{(k)}$ refine the very aspects of Schur functions 
that make them so fundamental and wide-reaching
and they are now called {\it $k$-Schur functions}.  

The role of $k$-Schur functions in affine Schubert
calculus emerged over a number of years.  The springboard was a
realization that the combinatorial backbone of $k$-Schur theory 
lies in the setting of the affine Weyl group.
The $k$-Schur functions are tied to Pieri rules, 
tableaux, Young's lattice, sieved $q$-binomial identities, and 
Cauchy identities that are naturally described in terms of 
posets of elements in $\tilde A^k$.  For example, 
$K_{\lambda,\mu}^{(k)}(1,1)$ is 
the number of reduced expressions for an element 
in $\tilde A^k$.
The combinatorial exploration fused into a geometric one when
the $k$-Schur functions were connected to the quantum cohomology 
of Grassmannians.
Quantum cohomology originated in string theory and symplectic geometry.
It has had a great impact on algebraic geometry and is intimately
tied to the {\it Gromov-Witten invariants}.  These invariants appear
in the study of subtle enumerative questions such as: how many degree 
$d$ plane curves of genus $g$ contain $r$ generic points?  
Lapointe and Morse \cite{[LMhecke]} showed that each
Gromov-Witten invariant for the quantum cohomology of Grassmannians 
exactly equals a $k$-Schur coefficient in the product of
{$k$-Schur functions} in $\Lambda$.  
A basis of {\it dual (or affine) $k$-Schur functions} was also 
introduced in \cite{[LMhecke]} and Lam proved \cite{Lam} that the 
Schubert bases for cohomology and homology of the affine 
Grassmannian $\Gr_{SL_{k+1}}$ are given by the
dual $k$-Schur functions and the $k$-Schur functions, respectively. 

Our motivation here is that the $k$-Schur functions $s_\lambda^{(k)}$
are parameterless and the $t$ is needed to connect with
theories outside of geometry.  Unfortunately, the characterizations  
for generic $t$ 
lack in mechanism for proofs.  We introduce a new family of functions that
reduce to $\{s_\lambda^{(k)}\}$ when $t=1$.
Our definition uses a combinatorial object called
affine Bruhat counter-tableaux (ABC's), whose
weight generating functions are the dual $k$-Schur functions \cite{DM}.
We associate a statistic (a non-negative integer) to each $ABC$
called the $k$-charge.  From this, we use the polynomials
\begin{equation}
\label{kkf}
K_{\lambda,\mu}^{(k)}(t)= \sum_{{\substack{\shape(A)=\core(\lambda) \\\ \weight(A)=\mu}}}
t^{k\text{-charge}(A)}
\end{equation}
to define a $t$-generalization of $s_\lambda^{(k)}$.
In particular,
we show that the matrix $(K^{(k)}_{\lambda,\mu}(t))_{\{\lambda,\mu
\in \mathcal P^k\}}$ is unitriangular and
taking the inverse of this matrix to be
$\tilde K^{(k)}$, a basis  for $\Lambda_t^k$ is given by
$$
s^{(k)}_\lambda[X;t] = 
\sum_{\mu} \tilde K_{\lambda,\mu}^{(k)}(t)\, H_\mu[X;t]\,,
$$
for all $\lambda$ with $\lambda_1\leq k$.
We prove that $s^{(k)}_\lambda[X;t]$ reduce to $k$-Schur
functions when $t=1$.  When $k=\infty$, these are
Schur functions, and thus \eqref{kkf} gives a new
description for the Kostka-Foulkes polynomials.
Naturally, we conjecture that these functions are
the $A_\lambda^{(k)}[X;t]$.

\section{Related work}

A refinement of the plactic monoid to a structure on $k$-tableaux
that can be applied to combinatorial problems involving $k$-Schur
functions is partially given in \cite{[LLMS2]} by a bijection
compatible with the RSK-bijection.  A deeper understanding of 
this intricate bijection is underway.  Towards this effort,
Lapointe and Pinto \cite{LaPi} have recently shown that a 
statistic on $k$-tableaux is compatible with the bijection.
There are now several statistics
(on $k$-tableaux, elements of the affine symmetric group,
and on $ABC$'s) whose charge generating functions are
the same.  The $ABC$'s can be used to find the image
of certain elements under this bijection and we are working
to put the $ABC$'s in a context that simplifies the bijection.

\section{Background}
\label{sec:tableauxIntro}
We identify each partition 
$\lambda=(\lambda_1,\ldots, \lambda_n)$ with its
Ferrers shape (having $\lambda_i$ lattice squares in the $i^{th}$ row, 
from the bottom to top).  For partitions $\lambda$ and $\mu$, we say $\lambda$ 
contains $\mu$, denoted $\mu \subseteq \lambda$, if $\lambda_i \geq \mu_i$.  A skew shape is a pair of partitions $\lambda$, $\mu$ such that $\mu \subseteq \lambda$, denoted $\lambda/\mu$.  

A \textit{semistandard tableau} $T$ is a filling of a Ferrers shape $\lambda$ with positive integers
that weakly decrease along rows and strictly increase up the columns.  The \textit{weight} of a semistandard
tableau is the composition $(\mu_i)_{i \in \mathbb{N}}$, where $\mu_i$ is the number of cells containing $i$.  For
a partition $\lambda$ and composition $\mu$, let $SSYT(\lambda,\mu)$ be the set of semistandard tableaux of 
shape $\lambda$ and weight $\mu$.

The \textit{hook-length} of a cell $(i,j)$ of any partition is the number of cells to the right of $(i,j)$ in row
$i$ plus the number of cells above $(i,j)$ in column $j$ plus 1.  A \textit{$p$-core} is a partition that 
does not contain any cell with hook-length $p$. 
The $p$-degree of a $p$-core $\lambda$,  $deg^p(\lambda)$,
is the number of cells in $\lambda$ whose hook-length is 
smaller than $p$.  Hereafter we work with a fixed integer $k > 0$ and all 
cores (resp. residues) are $k+1$-cores (resp. $k+1$-residues) and
$deg^{k+1}$ will simply be written as $deg$.  We let $\mathcal{P}^k$ denote the
set of all partitions $\lambda$ with $\lambda_1 \leq k$.  We also 
let $\mathcal{C}^{k+1}$ denote
the set of all $k+1$-cores.  We use a bijection given in \cite{[LMcore]}
$\mathfrak{c}: \mathcal{P}^k \rightarrow \mathcal{C}^{k+1}$.

For $n \geq 0$, an \textit{$n$-ribbon} $R$ is a skew diagram $\lambda/\mu$ consisting of 
$n$ rookwise connected cells such that there is no $2\times2$ shape contained in $R$.
We refer to the southeasternmost cell of a ribbon as its \textit{head}, and the northweasternmost 
cell of a ribbon as its \textit{tail}.

A \textit{ribbon tableau} $T$ of shape $\lambda/\mu$ is a chain of partitions 
\begin{displaymath}
\mu=\mu^0 \subset \mu^1 \subset \cdots \subset \mu^r = \lambda
\end{displaymath}
such that each $\mu^i/\mu^{i-1}$ is a tiling of ribbons filled with a positive integer.  A \textit{ribbon counter-tableau} $A$ of shape $\lambda/\mu$ is a ribbon tableau such that each skew shape $\mu^i/\mu^{i-1}$ is filled with the same positive integer $r-i+1$.  We set the cell $(i,j)$ of a ribbon counter-tableau to be the cell in row $i$, column $j$, where row one is the topmost row and column one is the leftmost column.  
For more on partitions and tableaux see \cite{Macbook}, \cite{Stanley}, \cite{FBergBook}.

\section{\texorpdfstring{Schubert representatives for $H^*(\Gr_{SL_{k+1}})$ and $H_*(\Gr_{SL_{k+1}})$}{Schubert representatives for H*(Gr(SL(k+1))) and H*(Gr(SL(k+1)))} }
\label{sec:strong k-pieri rule}


Despite the many characterizations for the Schubert representatives 
for the cohomology and homology of the infinite dimensional affine 
Grassmannian spaces for $SL_{k+1}$ (e.g. 
\cite{[LMcore],[LMhecke],[Lam],[LLMS],DM,AB}),
none have been shown to be the $t=1$ case of functions 
conjectured to give a positive Macdonald expansion \eqref{mackkostka}.
Our goal is to present functions with a $t$ parameter which
reduce to the $k$-Schur functions as formulated in \cite{DM} 
when $t=1$.  The formulation is given in terms of a combinatorial structure 
called $ABC$'s.

Recall that the strong (Bruhat) order on 
the affine Weyl group $\tilde{A}^{k}$ 
can be instead realized on $k+1$-cores 
by the covering relation:
\begin{displaymath}
  \rho \lessdot_B \gamma \Longleftrightarrow \rho \subseteq \gamma  
\text{ and } deg(\gamma) = deg(\rho)+1.
\end{displaymath}
An important fact about strong covers is useful in our study.

\begin{lemma}\cite{[LLMS]}
Let $\rho\lessdot_B\gamma$ be cores.  Then
\begin{enumerate}
\item Each connected component of $\rho/\gamma$ is a ribbon. 
\item The components are translates of each other and their heads
have the same residue. 
\end{enumerate}
\end{lemma}

A specific subset of ribbon counter-tableaux are those where each
ribbon is of height one.  An $ABC$ will be defined as such ribbon
counter-tableaux where the skew shapes are a certain strip defined
in terms of strong order.

\begin{defn}
\label{bsstrip}
For $0 < \ell \leq k$ and $k+1$-cores $\lambda$ and $\nu$, the skew shape $(k+\lambda_1,\lambda)/\nu$ 
is a \emph{bottom strong $(k-\ell)$-strip} if there is a saturated chain 
of cores
\begin{displaymath}
\nu = \nu^0 \lessdot_B \nu^1 \lessdot_B \cdots 
\lessdot_B \nu^{k-\ell}=(k+\lambda_1,\lambda)\,,
\end{displaymath}
where 
\begin{enumerate}
 \item $(k+\lambda_1,\lambda)/\nu$ is a horizontal strip
 \item The bottom rightmost cell of $\nu^i$ is also a cell in $\nu^i/\nu^{i-1}$, for $1 \leq i \leq k-\ell$.
\end{enumerate}
\end{defn}

It turns out that if a skew shape is a bottom strong strip then there is
a unique chain meeting the conditions described in Definition~\ref{bsstrip}.

\begin{exmpl}
The skew shape $(8,3)/(4,2)$ of 6-cores is a bottom strong 2-strip as there is the 
saturated chain
\begin{displaymath}
\text{\tiny\tableau[sbY]{&\cr&&&}} 
\lessdot_B
\text{\tiny\tableau[sbY]{&\cr&&&&&&}}
\lessdot_B
\text{\tiny\tableau[sbY]{&&\cr&&&&&&&}}\,.
\end{displaymath}
\end{exmpl}

\begin{exmpl}
The skew shape $(6,3,1,1)/(4,1,1,1)$ of 
4-cores is a bottom strong 1-strip as there is the saturated chain 
\begin{displaymath}
\text{\tiny\tableau[sbY]{ \cr \cr \cr & & & }} \lessdot_B \text{\tiny\tableau[sbY]{ \cr \cr & & \cr & & & & & }}\,.
\end{displaymath}
\end{exmpl}

\begin{exmpl}
There are 4 saturated chains of 4-cores in the strong order from $(3)$ to $(5,2,1)$,  
\begin{displaymath}
\text{\tiny\tableau[sbY]{ & & }} \lessdot_B \text{\tiny\tableau[sbY]{ \cr \cr & & }} \lessdot_B \text{\tiny\tableau[sbY]{ \cr \cr & & & }} \lessdot_B \text{\tiny\tableau[sbY]{ \cr & \cr & & & &}}\,,
\hspace{0.5in}
\text{\tiny\tableau[sbY]{ & & }} \lessdot_B \text{\tiny\tableau[sbY]{ \cr \cr & & }} \lessdot_B \text{\tiny\tableau[sbY]{ \cr & \cr & & }} \lessdot_B \text{\tiny\tableau[sbY]{ \cr & \cr & & & &}}\,,
\end{displaymath}
\begin{displaymath}
\text{\tiny\tableau[sbY]{ & & }} \lessdot_B \text{\tiny\tableau[sbY]{ \cr & & & }} \lessdot_B \text{\tiny\tableau[sbY]{ \cr \cr & & & }} \lessdot_B \text{\tiny\tableau[sbY]{ \cr & \cr & & & &}}\,,
\hspace{0.5in}
\text{\tiny\tableau[sbY]{ & & }} \lessdot_B \text{\tiny\tableau[sbY]{ \cr & & & }} \lessdot_B \text{\tiny\tableau[sbY]{ & \cr & & & & }} \lessdot_B \text{\tiny\tableau[sbY]{ \cr & \cr & & & &}}\,.
\end{displaymath}
Since none of these give a bottom strong strip, $(5,2,1)/(3)$ is not a bottom
strong strip.
\end{exmpl}

\begin{remark}
Bottom strong $(k-\ell)$-strips are a distinguished 
subset of strong strips in \cite{[LLMS]} that define the
Pieri rule for the cohomology of the affine Grassmannian.
\end{remark}

The iteration of bottom strong strips leads to the definition of an $ABC$.
First let us set some notation.
Given a ribbon counter-tableau $A$, let $A^{(x)}$ denote the subtableau made up of 
the rows of $A$ weakly higher than row $x$.
Let $A_{>i}$ denote the restriction of $A$ to letters strictly larger
than $i$ where empty cells in a skew are considered to contain
$\infty$.  With this in hand, we are now ready to define the $ABC$'s.

\begin{defn}
\label{def:abc}
For a composition $\alpha$ whose entries are not larger than $k$,
a skew ribbon counter-tableau $A$ is an \textit{affine Bruhat counter-tableau} 
(or $ABC$) of $k$-weight $\alpha$ if
$$(k+\lambda_1^{(i-1)},\lambda^{(i-1)})/\lambda^{(i)}\;\;\text{is a bottom strong $\alpha_i$-strip
\,for all $1\leq i\leq\ell(\alpha)$}\,,
$$
where
$\lambda^{(x)}=\shape(A^{(x)}_{>x})$.
We define the inner shape of $A$ to be $\lambda^{(\ell(\alpha))}$.
\end{defn}

The easiest method to construct an $ABC$ of $k$-weight $\alpha$ is iteratively,
from the empty shape $\lambda^{(0)}$, using Definition~\ref{bsstrip} to
successively add bottom strong strips that are
a tiling  of
$(k+\lambda^{(i-1)}_1,\lambda^{(i-1)})/\lambda^{(i)}$ with $\alpha_i$-ribbons
at each step.  


\begin{exmpl}
With $k=5$, we construct an $ABC$ of 5-weight $(3,3,1)$ by
\begin{displaymath}
\text{strong 3-strip}: 
\qquad
\text{\tiny\tableau[sbY]{&&}} 
\lessdot_B 
\text{\tiny\tableau[sbY]{&&&1}}
\lessdot_B 
\text{\tiny\tableau[sbY]{&&&1&1}}
\end{displaymath}
\begin{displaymath}
\text{strong 3-strip}: 
\qquad
\text{\tiny\tableau[sbY]{&\cr&&&}}
\lessdot_B
\text{\tiny\tableau[sbY]{&\cr&&&&{\color{blue}\bar 2}&{\color{blue}\bar 2}&{\color{blue}\bar 2}}}
\lessdot_B
\text{\tiny\tableau[sbY]{&&2\cr&&&&{\color{blue}\bar 2}&{\color{blue}\bar 2}&{\color{blue}\bar 2}&2}}
\end{displaymath}
\begin{displaymath}
\text{strong 1-strip}: 
\qquad
\text{\tiny\tableau[sbY]{&&\cr&&&}}
\lessdot_B
\text{\tiny\tableau[sbY]{3\cr&&\cr &&&&3}}
\lessdot_B
\text{\tiny\tableau[sbY]{3&3\cr&&\cr&&&&3&3}}
\lessdot_B
\text{\tiny\tableau[sbY]{3&3\cr&&\cr&&&&3&3&\color{red}\bar 3 & \color{red}\bar 3}}
\lessdot_B
\text{\tiny\tableau[sbY]{3&3\cr&&&3\cr&&&&3&3&\color{red}\bar 3&\color{red}\bar 3 & 3}}\,.
\end{displaymath}
The black letters are ribbons of size one, red letters make a ribbon of size two and blue letters make a ribbon of size 3 (or for those without color the ribbons are depicted with a bar).
This can be more compactly represented as
\begin{displaymath}
\text{\tiny\tableau[sbY]{3&3&2&1&1\cr&&&3&{\color{blue}\bar 2}&{\color{blue}\bar 2}&{\color{blue}\bar 2}&2 \cr &&&&3&3&{\color{red}\bar 3}&{\color{red}\bar 3}&3}}\,.
\end{displaymath}
\end{exmpl}

\begin{exmpl}
An example of an $ABC$ of $6$-weight $(4,4,2,1)$ with 
inner shape $(8,2,2,1) = \mathfrak{c}(6,2,2,1)$ is
\begin{displaymath}
\text{\tiny{ \tableau[sbY]{ &4&2&2&1&1 \cr & &3&3&3&{\color{blue}\bar 3}&{\color{blue}\bar 3}&{\color{blue}\bar 3}&2&2 \cr & &4&4&{\color{red}\bar 4}&{\color{red}\bar 4}&4&4&3&3&3&{\color{blue}\bar 3}&{\color{blue}\bar 3}&{\color{blue}\bar 3} \cr  & & & & & & & &4&4&{\color{red}\bar 4}&{\color{red}\bar 4}&4&4} }}\,. 
\end{displaymath}
\end{exmpl}

\begin{exmpl}\label{kstdABC}
Two examples of $ABC$'s of $k$-weight $(1,1,1,1,1,1,1) = (1^7)$ are
\begin{displaymath}
{\substack{{\text{\tiny{\tableau[sbY]{2&1&1 \cr 5&3&2&2 \cr  &4&{\color{red}\bar 3}&{\color{red}\bar 3}&3 \cr  &6&5&4&4 \cr & {\color{red}\bar 7}&{\color{red}\bar 7}&{\color{red}\bar 5}&{\color{red}\bar 5}&5 \cr  & & &6&{\color{red}\bar 6}&{\color{red}\bar 6} \cr & & &{\color{red}\bar 7}&{\color{red}\bar 7}&7}\,},}} \\\ {3\text{-weight }(1^7)}}}
\hspace{0.5in}
{\substack{{\text{\tiny{\tableau[sbY]{2&1&1&1&1&1&1 \cr 4&3&2&2&2&2&2&2 \cr 5&4&3&3&3&3&{\color{red}\bar 3}&{\color{red}\bar 3}&3 \cr 6&6&5&4&4&4&4&4&4 \cr  &7&7&5&5&5&5&{\color{red}\bar  5}&{\color{red}\bar 5}&5 \cr  & & &6&6&6&6&6&{\color{red}\bar 6}&{\color{red}\bar 6} \cr  & & & 7&7&{\color{red}\bar 7}&{\color{red}\bar 7}&7&7&7}\,}.}} \\\ {7\text{-weight }(1^7)}}}
\end{displaymath}
\end{exmpl}

The weight generating functions of the $ABC$'s turn out to be
the dual $k$-Schur functions.

\begin{theorem}
\cite{DM}
For any $\lambda \in \mathcal{C}^{k+1}$,
the dual $k$-Schur function can be defined by
\begin{displaymath}
\mathfrak{S}_{\lambda}^{(k)} = \sum_{A} x^{A}
\end{displaymath}
where the sum is over all affine Bruhat counter-tableaux of inner 
shape $\lambda$, and $x^{A} = x^{k\text{-}weight(A)}$.
\end{theorem}

These are symmetric functions, implying that
\begin{equation}
\mathfrak{S}_{\lambda}^{(k)} = \sum_{\mu: \mu_1\leq k}
K_{\lambda,\mu}^{(k)}\, m_\mu\,,
\end{equation}
where $K_{\lambda, \mu}^{(k)}$ is the number of affine Bruhat 
counter-tableaux of inner shape $\lambda$ and $k$-weight $\mu$. 
Then, using the Hall-inner product defined by
$$
\langle h_\lambda,m_\mu\rangle =\delta_{\lambda\mu}\,,
$$
we arrive at a characterization for $k$-Schur functions.
\begin{equation}
\label{kschurdef}
h_{\mu} = \sum_{\lambda} K_{\lambda,\mu}^{(k)}\,s_{\lambda}^{(k)}.
\end{equation}

\section{Kostka-Foulkes polynomials}
\label{sec:RWABC}

Our goal is to introduce polynomials $s_{\lambda}^{(k)}[X;t]$ 
that reduce to $s_{\lambda}^{(k)}[X]$ when $t=1$
using \eqref{mackkostka} as an inspiration.
Our approach is to introduce a statistic on 
$ABC$'s.
When $k=deg(\lambda)$, both $\mathfrak S_\lambda^{(k)}$ and
$s_\lambda^{(k)}$ are simply the Schur function $s_\lambda$.
One advantage of using $ABC$ combinatorics in the theory
of $k$-Schur functions is that known results concerning 
Schur functions can be reinterpreted in the $ABC$ framework
with $k$ large and this can shed light on the smaller
$k$ cases.  With this in mind, we consider a reformulation
for the Kostka-Foulkes polynomials in terms of 
$ABC$'s.  Our results enable us to give a 
characterization for symmetric polynomials in an extra
parameter $t$ that reduce to  
$\mathfrak S_\lambda^{(k)}$ and
$s_\lambda^{(k)}$ when $t=1$.

Let us start by recalling the Hall-Littlewood polynomials
$\{H_{\lambda}[X;t] \}_{\lambda}$.  
These are a basis for $\Lambda$ over the polynomial ring $\mathbb{Z}[t]$, 
which reduces to the homogeneous basis when the parameter $t=1$.
These often are denoted by $\{Q'_\lambda[X;t]\}$ in the
literature (\cite{Macbook}). Hall-Littlewood polynomials arise and can
be defined in various contexts such as the Hall Algebra, the character
theory of finite linear groups, projective and modular representations
of symmetric groups, and algebraic geometry.  We define them here via
a tableaux Schur expansion due to Lascoux and Sch\"utzenberger
\cite{LSfoulkes}.

The key notion is the \emph{charge} statistic on semistandard
tableaux. This is given by defining charge on words and then defining
the charge of a tableau to be the charge of its reading word.  For our
purposes, it is sufficient to define charge only on words whose
evaluation is a partition.
We begin by defining the charge of a word with weight $(1,1,\dots,1)$,
or a \emph{permutation}.  If $w$ is a permutation of length $n$, then
the charge of $w$ is given by $\sum_{i=1}^{n} c_i(w)$ where $c_1(w) =
0$ and $c_i(w)$ is defined recursively as
\begin{align*}
  c_i(w) &= c_{i-1}(w) + \chi\left( \text{$i$ appears to the right of
  $i-1$ in $w$} \right).
\end{align*}
Here we have used the notation that when $P$ is a proposition, $\chi(P)$ is
equal to $1$ if $P$ is true and $0$ if $P$ is false.

\begin{example} \label{Ex:ch1}
  The charge, $ch(3,5,1,4,2) = 0+ 1 + 1 + 2 + 2 = 6$.
\end{example}

We will now describe the decomposition of a word with partition
evaluation into \emph{charge subwords}, each of which are permutations.
The charge of a word will then be defined as the sum of the charge of
its charge subwords.  To find the first charge subword $w^{(1)}$ of a
word $w$, we begin at the \emph{right} of $w$ (i.e. at the last element
of $w$) and move leftward through the word, marking the first $1$ that
we see.  After marking a $1$, we continue to travel to the left, now
marking the first $2$ that we see.  If we reach the beginning of the
word, we loop back to the end.  We continue in this manner, marking
successively larger elements, until we have marked the largest letter
in $w$, at which point we stop.  The subword of $w$ consisting of the
marked elements (with relative order preserved) is the first charge
subword.  We then remove the marked elements from $w$ to obtain a word
$w'$.  The process continues iteratively, with the second charge
subword being the first charge subword of $w'$, and so on.
\begin{example} \label{Ex:ch2}
Given $w =(5,2,3,4,4,1,1,1,2,2,3)$, the first charge subword of $w$
are the bold elements in
  $(\mathbf{5}, \mathbf{2}, 3, 4, \mathbf{4}, 1, 1, \mathbf{1}, 2, 2,
  \mathbf{3})$.  If we remove the bold letters, the second
  charge subword is given by the bold elements in
$(\mathbf{3}, \mathbf{4}, 1, \mathbf{1}, 2, \mathbf{2})$.
It is now easy to see that the third and final charge subword is $(\mathbf{1}, \mathbf{2})$.  Thus we get that $ch(w) = ch(5,2,4,1,3) + ch(3,4,1,2) + ch(1,2)  = 8$.  Since $w$ is the reading word of the tableau
$
T=\text{\tiny{\tableau[scY]{ 5 \cr 2 & 3 & 4 & 4 \cr 1 & 1 & 1 & 1& 2 & 2 & 3}} }
$
we find that the $ch(T)=8$.
\end{example}

Equipped with the definition of charge, Hall-Littlewood polynomials
are then defined by
\begin{equation} \label{eq:Htos}
H_\mu[X;t] = \sum_\lambda
K_{\lambda,\mu}(t)  \,s_\lambda\,,
\end{equation}
where $K_{\lambda,\mu}(t)=K_{\lambda,\mu}(0,t)$ 
from \eqref{charge}.

Our first order of business to reformulate Kostka-Foulkes polynomials
is to describe the reading word of an $ABC$.
To do so, we first modify a given $ABC$ by lengthening the
row sizes. 

\begin{defn}
From a given $ABC$ $A$ of partition $k$-weight $\mu$, 
the extension of $A$, $ext(A)$, is the counter-tableau constructed 
from $A$ by adding $k$ cells with letter $i$ to each row $i$, where 
the first $\mu_i-s_i+r_i+1$ added cells form a ribbon
for $s_i$ the sum of the size of the ribbons filled with the 
letter $i$ in row $i$ and $r_i$ the number of such ribbons.  
\end{defn}

\begin{exmpl}\label{ssExtABC}
Consider the following extension of an $ABC$ with $5$-weight $(3,3,3,1)$.  
\begin{displaymath}
A = \text{\tiny\tableau[scY]{ &4&2&1&1 \cr &&4&3&{\color{blue}\bar 2}&{\color{blue}\bar 2}&{\color{blue}\bar 2}&2 \cr & & &{\color{red}\bar 4}&{\color{red}\bar 4}&4&{\color{red}\bar 3}&{\color{red}\bar 3}&3 \cr & & & & & &4&4&{\color{red}\bar 4}&{\color{red}\bar 4}&4}}
\Longrightarrow
ext(A) =  \text{\tiny\tableau[scY]{ &4&2&1&1&{\color{Maroon}\bar 1}&{\color{Maroon}\bar 1}&{\color{Maroon}\bar 1}&{\color{Maroon}\bar 1}&1 \cr  &&4&3&{\color{blue}\bar 2}&{\color{blue}\bar 2}&{\color{blue}\bar 2}&2&{\color{red}\bar 2}&{\color{red}\bar 2}&2&2&2 \cr  & & &{\color{red}\bar 4}&{\color{red}\bar 4}&4&{\color{red}\bar 3}&{\color{red}\bar 3}&3&{\color{blue}\bar 3}&{\color{blue}\bar 3}&{\color{blue}\bar 3}&3&3 \cr  & & & & & &4&4&{\color{red}\bar 4}&{\color{red}\bar 4}&4&4&4&4&4&4}}
\end{displaymath}
\end{exmpl}

\subsection{\texorpdfstring{Reading word of standard \textit{ABC}'s}{Reading word of standard ABC's}}
\label{sec:stdABC}

As with tableaux, we first define the reading word of a standard
$ABC$ (one of $k$-weight $1^n$) and use this to describe the
general reading word.  Standard $ABC$'s have a much more 
predictable structure than the general case.  Namely, a standard $ABC$ $A$ has only 
ribbons of size 1 or 2.  In fact, if a row $i$ in $A$ has an $i$-ribbon of size 2,
 then $\mu_i - s_i + r_i = 1$. Otherwise $\mu_i-s_i + r_i +1 = 2$.  
Thus, each row $i$ of $ext(A)$ has a unique $i$-ribbon of size 2.  

Our construction of the word of an $ABC$ $A$ considers only a subset of 
the cells in $ext(A)$.  Namely,
\begin{equation}
V_A = \{ (i,c_i)\in ext(A) :
\text{$(i,c_i)$ is any cell in a $i$-ribbon of row $i$ that is not its tail} \}.
\end{equation}
For standard $A$ of $k$-weight $1^n$, $V_A$ is simply a set of
$n$ ribbon heads; the one in each row $i$ of $ext(A)$ that 
contains $i$.  Using $V_A$, we define the reading word on this standard $ABC$ $A$.

\begin{defn} \label{stdWrdABC}
For a given $ABC$ $A$ of $k$-weight $1^n$, iteratively construct the 
reading word $w(A)$ by inserting letter $i$ directly right of letter 
$j$ where $j < i$ is the largest index 
such that $c_j < c_i$ and $(j,c_j) \in V_A$.
If there is no such $j$ then $i$ is placed at the beginning.
\end{defn}

\begin{exmpl}
Recall the $ABC$ from example \ref{kstdABC} of 3-weight $(1^7)$ is
\begin{displaymath}
A = \text{\tiny\tableau[scY]{2&1&1 \cr 5&3&2&2 \cr  &4&{\color{red}\bar 3}&{\color{red}\bar 3}&3 \cr  &6&5&4&4 \cr & {\color{red}\bar 7}&{\color{red}\bar 7}&{\color{red}\bar 5}&{\color{red}\bar 5}&5 \cr  & & &6&{\color{red}\bar 6}&{\color{red}\bar 6} \cr & & &{\color{red}\bar 7}&{\color{red}\bar 7}&7}}
\hspace{0.1in}
\Longrightarrow
\hspace{0.1in}
ext(A) = \text{\tiny\tableau[scY]{2&1&1&{\color{red}\bar 1}&{\color{red}\bar 1}&1 \cr 5&3&2&2&{\color{red}\bar 2}&{\color{red}\bar 2}&2 \cr  &4&{\color{red}\bar 3}&{\color{red}\bar 3}&3&3&3&3 \cr  &6&5&4&4&{\color{red}\bar 4}&{\color{red}\bar 4}&4 \cr & {\color{red}\bar 7}&{\color{red}\bar 7}&{\color{red}\bar 5}&{\color{red}\bar 5}&5&5&5&5 \cr  & & &6&{\color{red}\bar 6}&{\color{red}\bar 6}&6&6&6 \cr & & &{\color{red}\bar 7}&{\color{red}\bar 7}&7&7&7&7}}
\end{displaymath}
From $ext(A)$, we see that 
$
V_A = \{ (1,5),(2,6),(3,4),(4,7),(5,5),(6,6),(7,5)\}.
$
From $V_A$, we have the iterative construction of the reading word of $A$ as
$
(1) \rightarrow (1,2) \rightarrow (3,1,2) \rightarrow (3,1,2,4) \rightarrow (3,5,1,2,4) \rightarrow (3,5,6,1,2,4) \rightarrow (3,7,5,6,1,2,4).
$
This tells us that $w(A) = (3,7,5,6,1,2,4)$.
\end{exmpl}


\subsection{Reading words of ABC}
\label{sec:KFPolys}

Equipped with a method to obtain the reading word (permutation) of standard $ABC$'s,
we now define a way to construct a sequence of permutations from 
any $ABC$ with partition $k$-weight.

%
\begin{defn}\label{wordABC}
Let $A$ be an $ABC$ of partition $k$-weight $\mu$.  For $r=1,2,\ldots,\mu_1$, starting with $r=1$, we iteratively construct sets $E_A^r$ from $ext(A)$ as follows; put $(1,k + \mu_1 + 2 - r) \in E_A^r$, and let $(i,c_i) \in E_A^r$ if and only if $(i-1,c_{i-1}) \in E_A^r$ and
\begin{displaymath}
(c_{i-1} - c_i +k)~mod~(k+1) = min\{ (c_{i-1} - x+k)~mod~(k+1) | (i,x) \in V_A \setminus \displaystyle\cup_{p=1}^{r-1}E_A^{p} \}.
\end{displaymath}
\end{defn}
\begin{exmpl}\label{earABC}
Recall the $ABC$ $A$ from example \ref{ssExtABC} of $5$-weight $(3,3,3,1)$.
From its $ext(A)$, we see that 
\begin{displaymath}
V_A = \{ (1,7),(1,8),(1,9),(2,6),(2,7),(2,10),(3,8),(3,11),(3,12),(4,10)\}.
\end{displaymath}
We iteratively construct the sets $E_A^r$ for each $r=1,2,3$, using $ext(A)$ and $V_A$.  For $r=1$, begin by setting $E_A^1= \{(1,9)\}$.  Next, we see that $(2,7) \in E_A^1$ , because $(1,9) \in E_A^1$ and 
\begin{displaymath}
  1 = min \{ 2 = (9-6+5)~mod~6 , 1 = (9-7+5)~mod~6, 4 = (9-10+5)~mod~6 \}.
\end{displaymath}
So the next iteration gives us that $E_A^1 = \{(1,10),(2,7)\}$.  Next, we see that $(3,12) \in E_A^1$, because $(2,7) \in E_A^1$ and 
$
 0 = min \{ 4 = (7-8+5)~mod~6 , 1 = (7-11+5)~mod~6, 0 = (7-12+5)~mod~6 \}.
$
So the next iteration gives us that $E_A^1 = \{(1,10),(2,7),(3,12)\}$.  Finally since the $(4,10)$ is the only element in $V_A$ from the fourth row of $A$, then we see that
$
E_A^1 = \{ (1,9),(2,7),(3,12),(4,10)\}.
$

For $r=2$, to construct $E_A^2$, we begin by setting $E_A^2 = \{(1,8\}$, and repeat what we did to construct $E_A^1$, except this time we only consider elements from the set
$
V_A \setminus E_A^1 = \{(1,7),(1,8),(2,6),(2,10),(3,8),(3,11)\}.
$
This gives us 
$
E_A^2 = \{(1,8),(2,6),(3,11)\}.
$

Finally for $r=3$, to construct $E_A^3$, we begin by setting $E_A^3 = \{(1,7)\}$, and we only consider elements from the set
$
V_A \setminus (E_A^1 \cup E_A^2) = \{(1,7),(2,10),(3,8)\},
$
which immediately gives us
$
E_A^3 = \{(1,7),(2,10),(3,8)\}.
$
\end{exmpl}
Using each set $E_A^r$, we construct a sequence of reading word $w_r$ for $1 \leq r \leq \mu_1$.
\begin{defn}\label{ssReadWord}
Given an $ABC$ $A$ of partition $k$-weight $\mu$, for $1 \leq r \leq \mu_1$, the $r^{th}$ reading word of $A$, $w_r(A)$, is constructed using the same procedure in definition \ref{stdWrdABC}, where $V_A$ is replaced by
$E_A^r$.
\end{defn}

\begin{exmpl}\label{wordABCEx}
If we consider the $ABC$ $A$ from example \ref{ssExtABC}, then we know from example \ref{earABC} that
$
E_A^1 = \{ (1,9),(2,7),(3,12),(4,10)\}, ~
E_A^2 = \{(1,8),(2,6),(3,11)\},~
E_A^3 = \{(1,7),(2,10),(3,8)\}.
$
This tells us from definition \ref{ssReadWord} that $w_1(A) = (2,1,4,3)$, $w_2(A) = (2,1,3)$ and $w_3(A) = (3,1,2)$.
\end{exmpl}
For partitions $\lambda$, $\mu$ with $|\lambda| = |\mu| = n$, an $ABC$ $A$ is of $n$-weight $\mu$ and inner shape $\lambda$, has a charge statistic associated to it.

\begin{defn}
  Suppose $\lambda$ and $\mu$ are partitions with $|\lambda| = |\mu| = n$.  For any $ABC$ $A$ of $n$-weight $\mu$ and inner shape $\lambda$, the charge of $A$ is $ ch(A) = \sum_{r=1}^{\mu_1} ch(w_r(A))$.
\end{defn}

\begin{exmpl}\label{chABC}
The $ABC$ $A$ from example \ref{ssExtABC} has the reading words 
$
w_1(A) = (2,1,4,3),~w_2(A) = (2,1,3) \text{ and } w_3(A) = (3,1,2).
$
as described in example \ref{wordABC}.  Hence, we have that the charge of $A$ is 
$
ch(A) = ch((2,1,4,3)) + ch((2,1,3)) + ch((3,1,2)) = 2 + 1 + 2 = 5.
$
\end{exmpl}

There is a direct connection between reading words of semi-standard Young tableaux and a certain set of $ABC$'s.
\begin{theorem}
Suppose $\lambda$ and $\mu$ are partitions with $|\lambda| = |\mu| = n$. If the set
\begin{displaymath}
 ABC(\lambda,\mu) = \left\{A | ~ \text{$A$ is an $ABC$ of $n$-weight $\mu$ and inner shape $\mathfrak{c}(\lambda)$} \right\},
\end{displaymath}
then there is a bijection between the sets $ABC(\lambda,\mu)$ and $SSYT(\lambda,\mu)$, which is charge preserving.
\end{theorem}
From this theorem, we have the following corollary that gives the Kostka-Foulkes polynomials in the spirit of $ABC$'s.
\begin{cor}\label{KFABC}
  For partitions $\lambda$ and $\mu$, the Kostka-Foulkes polynomial
  \begin{displaymath}
    K_{\lambda,\mu}(t) = \sum_{A \in ABC(\lambda,\mu)} t^{ch(A)}.
  \end{displaymath}
\end{cor}

\section{\texorpdfstring{\textit{k}-charge and \textit{k}-Schur functions}{k-charge and k-Schur functions} }
\label{sec:StdkSchurExp}

We now look towards generalizing Corollary \ref{KFABC} by considering $ABC$'s of any partition $k$-weight.   
What is needed is an extra concept of an \textit{offset} of a given $ABC$.  Any $r$-ribbon of an $ABC$ is an \textit{offset} if there is a lower $r$-ribbon filled with the same letter as $R$ whose head has the same residue as the head of $R$.

\begin{defn}

For any $ABC$ $A$ of partition $k$-weight $\mu$, we set
\begin{displaymath}
  \text{off}~^k(A) = \sum_{\substack{\text{R: offset in }A}}\hspace{-0.15in}(\text{size}(R) - 1)
\end{displaymath}

\end{defn}

\begin{defn}
Let $A$ be an $ABC$ of partition $k$-weight $\mu$ and inner shape $\lambda$, and $w_1(A),\ldots,w_{\mu_1}(A)$ be the sequence of reading words that result from definition \ref{wordABC}.  Then, the $k$-charge of $A$
\begin{displaymath}
ch^k(A) = \sum_{r=1}^{\mu_1}ch(w_r(A)) - \textit{off}~^k(A)-\beta(A)
\end{displaymath}
where $\beta(A)$ is the number of cells in $\lambda$ whose hook-length exceeds $k$.
\end{defn}

\begin{exmpl}
Consider the $ABC$ of $3$-weight $1^5$ 
$
A = \text{\tiny\tableau[scY]{3&1&1 \cr  & 2 & {\color{red}\bar 2} & {\color{red}\bar 2}
\cr &4&3&3 \cr  &{\color{red}\bar 5}&{\color{red}\bar 5}&4&4 \cr & & & 5 &
{\color{red}\bar 5} & {\color{red}\bar 5} }}.
$
Here we see that $A$ has only one offset
$\text{\tiny\tableau[scY]{{\color{red}\bar 5}&{\color{red}\bar 5}}}$ in the second row from
the bottom.  The only reading word for this $A$ is $w_1(A) = (2,5,1,3,4)$.  So we get
$
ch^3(A) = ch((2,5,1,3,4)) - 1 - 1 = 5 - 1 - 1 = 3.
$
\end{exmpl}

\begin{defn}\label{genKFPoly}
For any $\lambda$, $\mu \in \mathcal{P}^k$,  we let
\begin{displaymath}
K_{\lambda,\mu}^{(k)}(t) = \sum_{\substack{\text{A: $ABC$ of $k$-weight $\mu$,} \\\ \text{inner shape $\mathfrak{c}(\lambda)$}}} \hspace{-0.35in} t^{ch^k(A)}.
\end{displaymath}
\end{defn}
Note that when $k \geq |\lambda|$, the polynomials $K_{\lambda,\mu}^{(k)}(t)$ are the Kostka-Foulkes polynomials of \ref{charge}.  Definition \ref{genKFPoly} generalizes the Kostka-Foulkes polynomials, and it also helps us to define a new set of symmetric functions with parameter $t$.  To see this, we only need the following claim.

\begin{lemma} \label{genKFMat}
The matrix
\begin{displaymath}
\left[ K_{\lambda,\mu}^{(k)}(t) \right]_{\{\lambda,\mu \in \mathcal{P}^k\}}
\end{displaymath}
is unitriangular.
\end{lemma}
Taking the inverse of the matrix in theorem \ref{genKFMat}, we form a basis for the subring $\Lambda_t^{(k)}$ of the ring $\Lambda$.

\begin{defn}\label{ktSchur}
For $\lambda \in \mathcal{P}^k$, the $k$-Schur function with parameter $t$ is
\begin{displaymath}
s_{\lambda}^{(k)}[X;t] = \sum_{\mu} \tilde{K}_{\lambda,\mu}^{(k)}(t)H_{\mu}[X;t],
\end{displaymath}
where $\tilde{K}_{\lambda,\mu}^{(k)}(t)$ are entries in the inverse of the matrix $\left[ K_{\lambda,\mu}^{(k)}(t) \right]_{\{\lambda,\mu \in \mathcal{P}^k\}}$.
\end{defn}

These new symmetric functions exhibit properties which connect them to the $k$-Schur and the Schur functions.
\begin{property}  As $\lambda$ ranges over partitions in $\mathcal{P}^k$, $s_{\lambda}^{(k)}[X;t]$ forms a basis for the subring $\Lambda_t^{(k)}$, $s_{\lambda}^{(k)}[X;1] = s_{\lambda}^{(k)}$, and $s_{\lambda}^{(\infty)}[X;1] = s_{\lambda}$.
\end{property}
Finally we make the following conjecture which ties the functions in Definition \ref{ktSchur} to those described in \cite{[LLM]}.

\begin{conj}
For $\mu \in \mathcal{P}^k$, 
$
s_{\mu}^{(k)}[X;t]  = A_{\mu}^{(k)}[X;t].
$
\end{conj}


\bibliographystyle{alpha}
\bibliography{fpsac13}

\begin{thebibliography}{LLMS12}

\bibitem[AB12]{AB}
S.~Assaf and S.~Billey.
\newblock Affine dual equivalence and $k$-{S}chur functions.
\newblock preprint, 2012.

\bibitem[Ber09]{FBergBook}
F.~Bergeron.
\newblock {\em {A}lgebraic {C}ombinatorics and {C}oinvariant {S}paces}.
\newblock A. K. Peters/CRC Press, 2009.

\bibitem[DM12]{DM}
A.~Dalal and J.~Morse.
\newblock The ${ABC}$'s of the affine {G}rassmannian and {H}all-{L}ittlewood
  polynomials.
\newblock {\em DMTCS Proceedings}, 2012.

\bibitem[GP92]{GP}
A.~M. Garsia and C.~Procesi.
\newblock {O}n certain graded ${S}_n$-modules and the $q$-{K}ostka polynomials.
\newblock {\em Adv. Math}, 87:82--138, 1992.

\bibitem[GR96]{GR}
A.~M. Garsia and J.~Remmel.
\newblock {P}lethystic formulas and positivity for $q,t$-{K}ostka coefficients.
\newblock {\em Mathematical essays in honor of Gian-Carlo Rota (Cambridge, MA},
  pages 245--262, 1996.

\bibitem[Gre55]{Green}
J.~A. Green.
\newblock {T}he characters of the finite general linear groups.
\newblock {\em Trans Amer Math Soc}, 80:442--407, 1955.

\bibitem[GT96]{GT}
A.~M. Garsia and G.~Tesler.
\newblock {P}lethystic formulas for {M}acdonald $q,t$-{K}ostka coefficients.
\newblock {\em {Adv Math}, {\bf}}, 123:144--222, 1996.

\bibitem[Hai01]{Haiman}
M.~Haiman.
\newblock {{H}ilbert schemes, polygraphs, and the {M}acdonald positivity
  conjecture}.
\newblock {\em J. Am Math. Soc.}, 14:941--1006, 2001.

\bibitem[KN96]{KiN}
A.~N. Kirillov and M.~Noumi.
\newblock $q$-difference raising operators for {M}acdonald polynomials and the
  integrality of transition coefficients.
\newblock {\em Algebraic methods and $q$-special functions (Montreal, QC},
  pages 227--243, 1996.

\bibitem[Kno97]{Knop}
F.~Knop.
\newblock {I}ntegrality of two variable {K}ostka functions.
\newblock {\em J. Reine Agnew}, 482:177--189, 1997.

\bibitem[Lam06]{[Lam]}
T.~Lam.
\newblock {A}ffine {S}tanley symmetric functions.
\newblock {\em Amer. J of Math}, 128(6):1553--1586, 2006.

\bibitem[Lam08]{Lam}
T.~Lam.
\newblock {S}chubert polynomials for the affine {G}rassmannian.
\newblock {\em J. Amer. Math Soc}, 21(1):259--281, 2008.

\bibitem[LLM03]{[LLM]}
L.~Lapointe, A.~Lascoux, and J.~Morse.
\newblock {T}ableau atoms and a new {M}acdonald positivity conjecture.
\newblock {\em Duke Math J}, 116(1):103--146, 2003.

\bibitem[LLMS10]{[LLMS]}
T.~Lam, L.~Lapointe, J.~Morse, and M.~Shimozono.
\newblock {A}ffine insertion and {P}ieri rules for the affine {G}rassmannian.
\newblock {\em Memoirs of the AMS}, 208(977), 2010.

\bibitem[LLMS12]{[LLMS2]}
T.~Lam, L.~Lapointe, J.~Morse, and M.~Shimozono.
\newblock {\em The poset of $k$-shapes and branching of $k$-{S}chur functions}.
\newblock to appear in Memoirs of the AMS, 2012.

\bibitem[LM05]{[LMcore]}
L.~Lapointe and J.~Morse.
\newblock {T}ableaux on $k+1$-cores, reduced words for affine permutations, and
  $k$-{S}chur function expansions.
\newblock {\em J Combin Theory Ser}, 112(1):44--81, 2005.

\bibitem[LM08]{[LMhecke]}
L.~Lapointe and J.~Morse.
\newblock {Q}uantum cohomology and the $k$-{S}chur basis.
\newblock {\em {Trans Amer Math Soc}, {\bf}}, 360(4):2021--2040, 2008.

\bibitem[LP]{LaPi}
L.~Lapointe and M.~E. Pinto.
\newblock Private communication.
\newblock Private communication with Authors.

\bibitem[LS78]{LSfoulkes}
A.~Lascoux and M.-P. Sch\"utzenberger.
\newblock {S}ur une conjecture de {H}.{O}. {F}oulkes.
\newblock {\em C.R. Acad. Sc. Paris}, 294:323--324, 1978.

\bibitem[Lus81]{Lu}
G.~Lusztig.
\newblock {S}ingularities, character formulas, and a $q$-analog of weight
  multiplicities.
\newblock {\em Analysis and topology on singular spaces, II, III (Luminy},
  101-102:208--229, 1981.

\bibitem[LV98]{LV}
L.~Lapointe and L.~Vinet.
\newblock A short proof of the integrality of the {M}acdonald $(q,t)$-{K}ostka
  coefficients.
\newblock {\em Duke Math J}, 91:205--214, 1998.

\bibitem[Mac88]{[M2]}
I.~G. Macdonald.
\newblock {\em A new class of symmetric functions}.
\newblock S\'eminaire Lotharingien de Combinatoire, {\bf B20a} 41pp, 1988.

\bibitem[Mac95]{Macbook}
I.~G. Macdonald.
\newblock {\em {S}ymmetric functions and {H}all polynomials}.
\newblock Clarendon Press, Oxford, 2nd edition, 1995.

\bibitem[NY97]{NY:1997}
A.~Nakayashiki and Y.~Yamada.
\newblock {K}ostka polynomials and energy functions in solvable lattice models.
\newblock {\em Selecta Math}, 3(4):547--599, 1997.

\bibitem[Sah96]{S}
S.~Sahi.
\newblock {I}nterpolation, integrality, and a generalization of {M}acdonald
  polynomials.
\newblock {\em Internat Math. Res. Notices}, pages 457--471, 1996.

\bibitem[Sta99]{Stanley}
R.~Stanley.
\newblock {\em {E}numerative {C}ombinatorics}.
\newblock Cambridge, Vol 2, 1999.

\end{thebibliography}
\label{sec:biblio}

\end{document}